## Paper information

| | |
|---|---|
| Authors | Bin Liu[◊, †], Thomas Brinsmead[†], Stefan Westerlund[‡], Robert Davy[‡] |
| Affiliations (optional) | [◊]Network Planning Division, Transgrid, Sydney 2000, Australia |
| | [†]Energy Centre, CSIRO, Mayfield West 2304, Australia[1] |
| | [‡]Information Management & Technology, CSIRO, Black Mountain 2601, Australia |
| Email address | Bin Liu (Email: bin.liu@transgrid.com.au, eeliubin@ieee.org) |

## Summary[2]

Worldwide renewable generation capacity has been increasing rapidly over the past decade. Unlike conventional generators, renewable energy sources such as wind and solar farms are less controllable, leading to increasing uncertainty in power system operation. It is important that a power system operator can take measures to ensure power system reliability and security even under such circumstances.

One approach to assessing the security of a power system operating point with respect to uncontrollable changes in renewable generation quantities is by calculation of the dispatchable region for renewables (DRR). This is a *region* within the space of possible renewable generation quantities that a power system operator can satisfactorily manage by dispatching available controllable resources. Geometrically, the DRR is a polyhedron with boundaries determined by binding operational limits on generation and transmission. It can be used: a) to evaluate how close a power system operating point is to a boundary of the secure operational region, and b) to identify the uncontrollable generation ramping event that is the shortest distance to the boundary.

The DRR can be calculated by solving a series of *max-min* problems based on constraints representing the power system's operational limits, which include the output and ramping limits of controllable generators and power flow limits on transmission lines. One common approach to solving each *max-min* optimisation is through mixed-integer linear programming (MILP) reformulation. However, the resulting MILP can be computationally burdensome for a large-scale power system with many renewable farms (non-dispatchable generation nodes). An alternative approach is called the iteration-based linear programming (IBLP) algorithm, where the *max-min* problem is solved by alternately solving *max* and the *min* problems with numerous candidate solution starting points. However, a large number of starting points may also be computationally expensive.

This paper investigates whether advanced computation, including both high-performance computing (HPC) and parallel computing techniques (PCT), can improve computational performance in calculating DRR. Specifically, standard power system test networks, namely the PJM-5, IEEE-39 and IEEE-118 bus systems, will be used as the basis for building dispatch models with various numbers of renewable farms using PowreModels.jl. The calculation of DRR based on both MILP and IBLP approaches will be tested on two hardware platforms, one with and one without advanced computational techniques.

---

[1] Part of the work was completed when the author was with Energy Centre, CSIRO, Mayfield West 2304, Australia.

[2] Disclaimer: The views expressed herein are not necessarily the views of any organisations affiliated to the authors, and the affiliated organisations of all authors do not accept responsibility for any information or advice contained herein.

## Acronyms

| | | | |
|---|---|---|---|
| **BESS** | Battery energy storage system | **KKT** | Karush-Kuhn-Tucker |
| **CPT-01/02** | Computation platform-01/02 | **LP** | Linear programming |
| **DC-OPF** | Direct-current optimal power flow | **MILP** | Mixed-integer linear programming |
| **DRR** | Dispatchable region for renewables | **MIP** | Mixed-integer linear |
| **FCAS** | Frequency control ancillary services | **PCT** | Parallel computation techniques |
| **HPC** | High-performance computing | **PSD** | Power system dispatch model |
| **IBLP** | Iteration-based linear programming | **SOS1** | Special order sets of type 1 |
| **IBLP-ST/MT** | IBLP with single/multiple threads | | |

## 1. Introduction

Worldwide renewable generation capacity has been increasing rapidly over the past decade. Compared to conventional synchronous generators, large-scale renewable farms are less controllable, increasing uncertainty in power system operations. Moreover, increasing penetration of distributed resources (DERs) makes it more difficult to accurately predict loads at the transmission-distribution interface. To cope with such uncertainties, it is critical that a power system operator can take actions as necessary to ensure power system reliability and security, especially when generation outputs from renewable farms are close to or exceed the operational limits of the rest of the power system.

One approach to assessing the security of a power system operating point with respect to possible uncontrollable changes in renewable generation quantities is by calculation of the dispatchable region for renewables (DRR) [1,2]. This is the *region* of permissible (uncontrollable) renewable generation quantities that a power system operator can manage by dispatching existing controllable resources. Such controllable resources can be conventional synchronous generators, battery energy storage systems (BESSs) or controllable demand. However, the capability of a power system to cope with uncertainties in renewable generation output is constrained by several factors, including both the ramping capability and output capacity limits of controllable generators, and also the power flow limits on transmission lines.

When a linear transmission network model, such as that in direct-current optimal power flow (DC-OPF), is used to calculate generation dispatch constraints, the DRR can be geometrically represented as a polyhedron, with boundaries determined by binding power system operational constraints. The relevance of a DRR is that if actual renewable generation outputs fall within it, the power system operator can always find a dispatch strategy for the remaining controllable generators to meet demand while ensuring no operational constraints are violated. Therefore, the DRR can be used: a) to monitor how close a given operating point for renewable generation is to the boundary of a secure operational region, and b) to identify the renewable generation ramping event that is the closest to that boundary.

However, calculating DRR can be challenging, particularly when an *affine* control policy, which pre-defines how controllable generators will be controlled to ensure supply-demand balance when renewable generations deviate from the forecasted values, is not already known and available [3]. In the absence of an *affine* control policy, it is instead possible to identify a series of successively smaller closed sets that contain the DRR by solving a series of *max-min* problems. Such a calculation procedure will be discussed in this paper in detail.

For solving the relevant *max-min* problem, either a mixed-integer linear programming (MILP)-based approach (MILP approach) or an iteration-based linear programming (IBLP)-based approach (IBLP approach) can be used, as discussed in [1,2,4]. These studies show that the computational burden of the MILP approach can be high for large-scale systems with many renewable farms, even with the special ordered sets of type 1 (SOS1) techniques supported by some commercial solvers like Cplex [5]. In contrast, the IBLP approach, although unable to guarantee a global optimum, can provide solutions of high quality with reasonable computational time requirements, provided the algorithm is given a sufficient number of carefully selected starting points [2,6].

This paper does not propose or explore novel methods for solving the *max-min* problem. Instead, it experimentally investigates whether advanced computation techniques, including both high-performance computing (HPC) and parallel computation techniques (PCT), can improve the computational performance of calculating DRR, and provides a perturbation-based approach to identify the binding power system operational constraints corresponding to DRR boundaries. Note that: 1) The MILP approach, due to its complexity, is computationally demanding. Whether and when HPC can be beneficial needs further investigation; and 2) In the IBLP approach, several independent linear programming (LP) problems must be solved for each (iteration of the) *max-min* problems and, however, given multiple CPU threads, the LP problems can be solved in parallel.

The remainder of the paper is organised as follows. The canonical power system dispatch model (PSD) is described in Section 2. Following this in Section 3, we will present the definition of DRR, revisit approaches to calculating DRR, describe the calculation workflow for a Julia Programming environment, present the perturbation-based approach to identify the binding power system operational constraints for each DRR boundary, and discuss how advanced computational techniques can be used. Case studies based on the standard networks: the PJM-5 bus, IEEE-39 bus, and IEEE-118 bus systems, will be presented in Section 4 and concluding remarks are provided in Section 5.

## 2. Canonical Power System Dispatch Model

In this paper, the DRR will be defined for a canonical PSD, formulated as in [1] as

$$\min_{p_i} \sum_i f(p_i) \tag{1a}$$

$$s.t. \quad \sum_i p_i + \sum_i w_i = \sum_i d_i \tag{1b}$$

$$-L_l^+ \leq \sum_i \pi_{il}(p_i + w_i - d_i) \leq L_l^+ \quad \forall l \tag{1c}$$

$$\max\{p_i^* - R_i^-, p_i^-\} \leq p_i \leq \min\{p_i^* + R_i^+, p_i^+\} \quad \forall i \tag{1d}$$

where $w_i$ is the actual generation (active power) realised by renewable farm $i$; $p_i$ is pre-dispatched controllable generation output from the $i^{th}$ flexible unit while $p_i$ is the updated dispatched power after observing the $w_i$; $d_i$ is the active power demand at bus $i$; $f(p_i)$ represents the objective function related to $p_i$; $R_i^-/R_i^+, p_i^-/p_i^+$ are downward/upward ramping rates and lower/upper generation limits of unit $i$, respectively; $\pi_{il}$ is the sensitivity of power injection at bus $i$ to active power flow through line $l$ and $L_l^+$ is the power flow limit on transmission line $l$.

Note that flexible units here (with power $p_i$) could be either conventional synchronous generators or BESSs that have the capability to adjust their generated power (for example, when $w_i$ deviates from the forecasted value). However, an inflexible unit can alternatively be represented as a *flexible* unit with zero ramping capability.

In the PSD model, the objective function (1a) is subject to the power balance constraint (1b), power flow constraints on transmission lines (1c), and ramping capability constraints for each flexible generator in (1d). Fortunately, other constraints, e.g., voltage angle difference constraints, which by default are considered in the DC-OPF model in PowerModels.jl, can also be conveniently represented in the PSD formulation.

It is noteworthy that the PSD model may be different from those employed for operating a real power system, particularly in an electricity market environment. For example, the Australian Energy Market Operator (AEMO) co-optimises both the energy and Frequency Control Ancillary Services markets in the National Electricity Market (NEM) to minimise the total operational cost while considering network constraints [7]. In such a case, (1) can be revised to take into account both the energy and FCAS offers in the objective function, and constraints should cover both the coupled energy/FCAS constraints (a variant of (1d)), and the FCAS demands. Nevertheless, the approach to calculating DRR can be conveniently extended to and applied in a real power system if the required energy and FCAS data are available.

For the convenience of later discussion, we here express the PSD model in a more compact form. Define vector variables $p, w, v$ and $l$ representing, respectively: the active powers from flexible generators, active powers from renewable farms, nodal voltage angles, and active power flows in transmission lines. The PSD model can be expressed in compact form as

$$\min_{p,l,v} f(p) \tag{2a}$$
$$s.t. \quad Aw + Bp + Cl + Dv = b \tag{2b}$$
$$Ew + Fp + Gl + Jv \leq d \tag{2c}$$

where $A, B, C, D, E, F, G$ and $J$ are known matrices with appropriate dimensions, and similarly $b$ and $d$ are known vectors.

## 3. Definition and calculation of DRR

Considering (2), it is observed that the constraints (2b)-(2c) define a feasible region for $(w, p, l, v)$, which is a polyhedron. However, since $w$ is less controllable (or uncontrollable), a system operator may be interested in the *space*, within which $w$ falls, there is always a controllable $(p, l, v)$ such that (2b)-(2c) are not violated. Such a *space* for $w$ (a polyhedral subset of the projected subspace of $(w, p, l, v)$) defines the DRR to be studied. In other words, the objective function in the PSD model is *not* relevant to the calculation of DRR. Further, a *max-min* optimisation problem can be formulated, where an initial set containing the DRR for $w$ will be iteratively updated until its optimal objective value equals 0. The mathematical definition of the DRR and its calculation procedure will be discussed in this section.

### 3.1 Revisiting approaches to calculating DRR

#### 3.1.1 Concept and mathematical model of DRR

The DRR, as discussed in [1] and [2], can be defined as a maximum polyhedron, where for any realisation of *w* within the polyhedron, there is always a solution $(p, l, v)$ such that (2b)-(2c) are not violated. In other words, we have

$$\mathrm{DRR} = \{w | \exists (p, l, v) \ s.t. \ Aw + Bp + Cl + Dv = b, Ew + Fp + Gl + Jv \leq d\} \tag{3}$$

Note that the DRR is a function of known parameters, including the pre-dispatch set-point $p^*$ that may be updated regularly. Moreover, DRR is the projection of the feasible region for $(w, p, l, v)$ on the subspace only containing $w$, and thus- as a polyhedron- it can alternatively be expressed as $\mathrm{DRR} = \{w | Hw \leq f\}$ for some $H$ and $f$. However, as discussed in [1,2], algorithms to directly identify the *projected space* can be computationally intractable for large systems. Alternatively, starting from an initial polyhedron, say $W_0$, that is sufficiently large so that it encloses the DRR, the DRR can be subsequently identified (in principle, exactly) by iteratively adding linear cuts after identifying one or multiple points within $W_k$ (the updated $W$ starting from $W_0$ in the $k^{th}$ iteration), that fall outside the DRR. Whether or not $W_k$ contains points outside the DRR can be checked by solving the *max-min* problem (4) following after introducing non-negative slack variables $t^+, t^-$ and $s$ to equations (2b)-(2c) as follows.

$$F_k = \max_{w \in W_k} \min_{p,l,v,t^+,t^-,s} 1^T s + 1^T (t^+ + t^-) \tag{4a}$$
$$s.t. \quad Aw + Bp + Cl + Dv + t^+ - t^- = b \quad (\alpha) \tag{4b}$$
$$Ew + Fp + Gl + Jv - s \leq d \quad (\delta) \tag{4c}$$
$$s \geq 0, t^- \geq 0, t^+ \geq 0 \tag{4d}$$

where $\alpha$ and $\delta$ are dual variables for corresponding constraints.

Observe that $F_k > 0$ implies that $W_k$ contains points that fall outside of DRR and $W_k$ needs to be further updated. Otherwise, if after several iterations, we have $F_k = 0$, then $W_k$ is the target DRR. It is noteworthy that whenever there exists at least one renewable operating point $w$ for which a dispatch strategy can be found, i.e., whenever DRR is non-empty, it is expected that $F_k = 0$ for some finite number $k$ of iterations.

We now discuss how to solve the *max-min* problem (4) and how to update from $W_k$ to $W_{k+1}$.

### 3.1.2 Reformulating and solving the *max-min* problem

As (4) is a strongly non-convex problem, which is difficult to solve directly, one common practice is to first formulate the dual of the *min* (sub)problem and solve the resulting *max* problem. Exploiting duality, an equivalent formulation of (4) is given below, dropping the subscript $k$ for convenience. More details can be found in Appendix 6.

$$F = \max_{\alpha,\delta} \alpha^T b - \delta^T d + \max_{w}(-\alpha^T A + \delta^T E)w \tag{5a}$$
$$s.t. \ \ 0 \leq \delta \leq 1, -1 \leq \alpha \leq 1 \tag{5b}$$
$$\alpha^T B - \delta^T F = 0, \alpha^T C - \delta^T G = 0, \alpha^T D - \delta^T J = 0 \tag{5c}$$
$$Hw \leq f \tag{5d}$$

Note that (5) is still a non-convex problem, and it now has bilinear terms in both the objective function and the linear constraints. Two typical approaches to solving (5) are revisited following.

**The mixed-linear integer programming (MILP) approach.** In the MILP approach, the optimisation problem (5) with the bilinear objective function is first reformulated as a MILP problem. The following steps exploit Karush-Kuhn-Tucker (KKT) conditions for extrema.

Observe that if $\alpha$ and $\delta$ are fixed, the inner *max* problem in (5), i.e., $\max_{w}(-\alpha^T A + \delta^T E)w | s.t. Hw \leq f\}$, is an LP problem. The KKT conditions for this LP problem are

$$\frac{\partial L}{\partial w} = A^T\alpha - E^T\delta + H^T\theta = 0 \tag{6a}$$
$$0 \leq \theta \perp f - Hw \geq 0 \Leftrightarrow \theta^T Hw = \theta^T f \tag{6b}$$

where $L = (\alpha^T A - \delta^T E)w + \theta^T(Hw - f)$ is the Lagrangian function of the LP problem.

Rearranging (6a) for $\theta^T H$ and substituting in (6b), we conclude that $(-\alpha^T A + \delta^T E)w = \theta^T Hw = \theta^T f$. Consequently, (5) can be reformulated as the following optimisation problem with complementary constraints.

$$\max_{\alpha,\delta,\theta,w} \alpha^T b - \delta^T d + \theta^T f \tag{7a}$$
$$s.t. \ \ (5b) - (5c), (6a) \tag{7b}$$
$$0 \leq \theta \perp (f - Hw) \geq 0 \tag{7c}$$

Note that (7c) is equivalent to $\theta_i(f_i - H_{i,:}w) = 0$ for each $i$. These conditions can be further reformulated as a set of mixed-integer linear (MIP) constraints. After introducing binary variables $\lambda_i$ and a large constant $M$, reformulating (7c) leads to

$$0 \leq \theta_i \leq (1 - \lambda_i)M \tag{8a}$$
$$0 \leq f_i - H_{i,:}w \leq \lambda_i M \tag{8b}$$
$$\lambda_i \in \{0,1\} \tag{8c}$$

Note that the computational burden of solving the MILP problem can be reduced if redundant constraints characterising $W_k$ are removed each time $W_k$ is updated.

**The iterative linear programming approach (IBLP).** In this approach, (5) will be solved by selecting several initial fixed guesses, say $w_j$ for the $j^{th}$ guess (that is, $w_{kj}$ as the $j^{th}$ guess in the $k^{th}$ iteration), and alternately solving the *max* (sub)problem for $(\alpha,\delta)$ (similarly, these are $(\alpha_{kj},\delta_{kj})$ as the $j^{th}$ guesses in the $k^{th}$ iteration) by fixing $w$ and the *max* (sub)problem for $w$ by fixing $(\alpha,\delta)$, iterating until the optimal objective values of the two optimisation problems are sufficiently close to each other. However, the solution quality depends strongly on both the number and the quality of the initial estimates $w_j$. Possible approaches to generating effective $w_j$ have been discussed for $W$ in polyhedron form in [2] and in ellipsoidal form in [6]. In this paper, candidate $w_j$ for solving (5) are generated as follows.

1. Generating set $S_1$: points are randomly generated within the initialised $W_0$, which means $S_1$ will be the same for all outer iterations (that is, $S_1$ is independent of $k$)[3].
2. Generating scenario set $S_2$: points are generated by projecting $\overline{w}$, the forecasted or expected renewable generation vector, to the boundaries of $W$ (that is, $W_k$) at the beginning of each iteration[4]. Specifically, for each halfplane in $W$, say $H_{i,:}w \leq f_i$, where $H_{i,:}$ and $f_i$ are the $i^{th}$ row of $H$ and the $i^{th}$ element of $f$, respectively, the scenario point is generated as $w = \overline{w} - \frac{H_{i,:}^T \overline{w} - f_i}{H_{i,:}^T H_{i,:}} H_{i,:}$.
3. Generating scenario set $S$: the set is generated as $S = S_1 \cup S_2$, which will be used for solving (5).

Note that as $W$ is updated iteratively, $S_2$, and hence $S$, will also change with each iteration $k$.

### 3.1.3 Generating linear cuts for updating $W$

With optimal solutions of $\alpha^*, \delta^*$ and $w^*$ reported by either IBLP or MILP approach, and if $F_k^* > 0$, the following linear cuts will be generated to update $W_{k+1}$.

$$(\alpha^*)^T b - (\delta^*)^T d + [-(\alpha^*)^T A + (\delta^*)^T E]w \leq 0 \quad (9)$$

Note that for the MILP approach, only one linear cut will be added to $W$ in each iteration $k$, since only one worst $w$ can be identified. However, for the IBLP approach, several linear cuts can be added in each iteration since there might be more than one initial guess of $w$ (that is, across several $j$) that leads to a positive optimal objective value. By adding more linear cuts in each iteration, the IBLP approach is expected to take fewer $k$ iterations than the MILP approach.

## 3.2 Workflow for DRR calculation

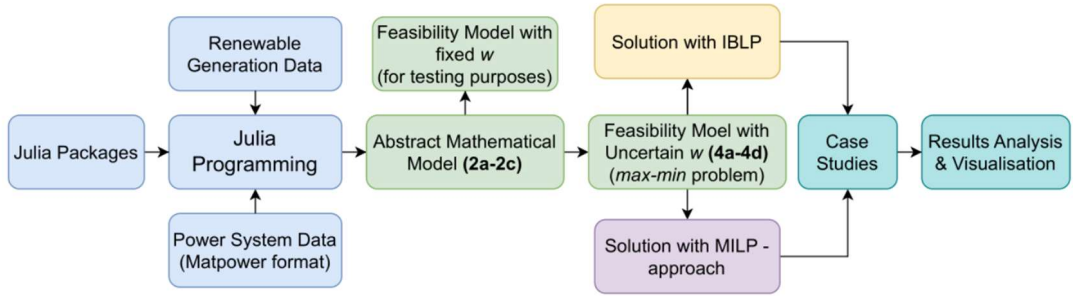

*Figure 1: The schematic overview of the DRR calculation.*

We have realised the calculation of DRR in Julia, which is a programming language especially designed for scientific computing [8] and supported by packages including PowerModels.jl [9], JuMP.jl [10], LinearAlgebra.jl, MathOptInterface.jl [11], JLD2.jl, Gurobi.jl and Plots.jl etc. Power system models were provided in Matpower format [12]. The workflow we implemented for calculating DRR appears in Fig.1.

Several remarks on the workflow are given below.

1. The original Matpower data file is modified to include renewable farms, each of which is represented as a single generator.
2. With the assistance of MathOptInterface.jl, the abstract model (2) is extracted from the DC-OPF formulation defined in PowerModels.jl after defining each generator's ramping capability as 10% of its maximum available active power.

---

[3] Alternatively, $S_1$ could be updated at the beginning of each iteration $k$. However, for consistency in comparing simulations, using single-thread or multiple threads, or across alternative computational platforms, $S_1$ is fixed in all iterations over $k$.

[4] $W$ should be interpreted as the updated $W$, i.e., $W_k$, at the $k^{th}$ iteration.

3. The abstract model is tested against the initially defined DC-OPF model in the block: *Feasibility Model with Fixed $w$ (for testing purposes)*.
4. For *Results Analysis & Visualisation*, DRR is only provided when its dimension count is two or less, while the figure showing how the *max-min* problem objective value changes along with the iteration process will be presented for all studied cases.

### 3.3  Identifying binding power system constraints for DRR boundaries

As DRR boundaries are identified by minimising the constraint violations via (5), each of the identified boundaries is connected to one or more *binding constraint(s)*. Although such binding constraints can be identified through the DRR calculation process, mapping such information in the compact model, which is derived via Julia packages, back to the physical PSD model, which is defined in Matpower format and formulated via PowerModels.jl, can be complicated. Instead, we propose here to identify the binding constraints by a perturbation-based approach post the DRR calculation, as we explain next.

With the DRR expressed as $\{w|Hw \leq f\}$, the terminal points[5] for each of its boundary hyperplanes can be calculated by solving the optimisation problem: $\max_w / \min_w \{\beta^T w | Hw \leq f, H_{i,:}w = f_i\}$ with $\beta$ being a properly defined constant vector, which by default can be set as $\beta = \mathbf{1}$. Assuming the two identified terminal points for the $i^{th}$ constraint are expressed as $w^{u,i}$ and $w^{l,i}$, the binding constraints can be identified by the following steps.

1. Building the *feasibility model* for (1), where the (1c) and (1d) are slacked by introducing non-negative variables and the objective function (1a) will be replaced by minimising the sum of all slack variables.
2. Fixing $w$ at $\frac{w^{u,i}+w^{l,i}}{2} + \lambda\left(\frac{w^{u,i}+w^{l,i}}{2} - \overline{w}\right)$ and solving the *feasibility model*, where $\lambda$ is a scaling factor to locate a perturbation point that is outside of DRR and that is sufficiently small, in order to identify the binding constraints. In this paper, we have $\lambda = 0.01$.
3. If the optimal objective value for solving the *feasibility model* is positive, the binding constraints can be identified as those constraints for which optimal values of slack variables are reported as positive. Otherwise, it is implied that the binding constraints are determined by $W_0$ and can be identified by comparing $\frac{w^{u,i}+w^{l,i}}{2}$ and the boundaries of $W_0$.

Note that different binding constraints may lead to the same DRR boundary, as we explain with an illustrative example. Consider a system with 2 buses, where Generator 1 (dispatched at 70 MW with a ramping capability of $\pm 10$ MW) is connected to bus 1, Generator 2 (dispatched at 30 MW with a ramping capability of $\pm 5$ MW), a constant load of 150 MW and a renewable farm (generating at 50 MW) are connected to bus 2, the line connecting bus 1 and bus 2 is thermally limited by 70 MW. Obviously, the DRR for this renewable farm is $\{w|45 \leq w \leq 65\}$. For the upper bound, i.e., $w \leq 65$, the binding constraints are ramp-down capabilities of both Generator 1 and Generator 2. By contrast, for the lower bound, i.e., $w \geq 45$, the binding constraints are ramp-up capability of Generator 2 and the thermal limit of the connecting line. For such a case, there should be multiple optimal solutions for the *feasibility model*.

However, only a single binding constraint may be identified since the solver generally reports only one optimal solution. An approach to further refine the solution is by properly defining the penalty coefficients when building the *feasibility model* for (1), where constraint violations can be allocated with various weights (*prices*) to represents the severities of the violations and/or the system operator's experience based preferences.

### 3.4  The applications of advanced computing techniques

As discussed in Section 3.1, advanced computational techniques, including both HPC and PCT, can be applied to estimate DRR by solving the series of *max-min* problems in order to identify valid linear cut(s).

---

[5] The *terminal points* of a boundary halfplane refer to the extrema of the DRR (a polyhedron) that also belong to this boundary halfplane.

In the MILP approach, a series of MILP problems, with possibly numerous binary variables, must be solved, which may require significant computational resources. This suggests that HPC could accelerate the computational performance of the MILP approach. In contrast with the IBLP approach, numerous LP problems (more than one for each initial guess $w_j$) must be solved in each iteration. Although LP problems are generally much easier to solve than MILP problems, the number of LP problems to be solved can be larger for the IBLP approach. However, because the LP problems are independent for each $j$, this suggests that PCT may improve the computational performance of the IBLP approach.

In order to investigate improvements due to HPC and PCT, three power system networks, including the PJM-5 bus, IEEE-39 bus, and IEEE-118 bus, systems were each studied both on: a) CPT-01: a laptop with 11th Gen Intel(R) Core(TM) i7-1185G7@3.0 GHz 4-Core CPU and 16 GB RAM, and on b) CPT-02: a remote server with 2 × Intel Xeon Platinum 8280L@2.7 GHz 28-Core CPUs and 3 TB RAM. For this paper, all MILP and LP problems were solved by Gurobi [13].

When investigating the performance improvement due to PCT, 4 threads and 50 threads are enabled for CPT-01 and CPT-02, respectively.

## 4. Case Study

### 4.1 Case setup

Generator and network data of the original PJM 5-bus system can be found in [14] and is also available in attached data sets in PowerModels.jl. For this paper, two renewable farms are added, one each to bus 4 and bus 5. Some parameter values for the modified PJM 5-bus system appear below in Fig.2.

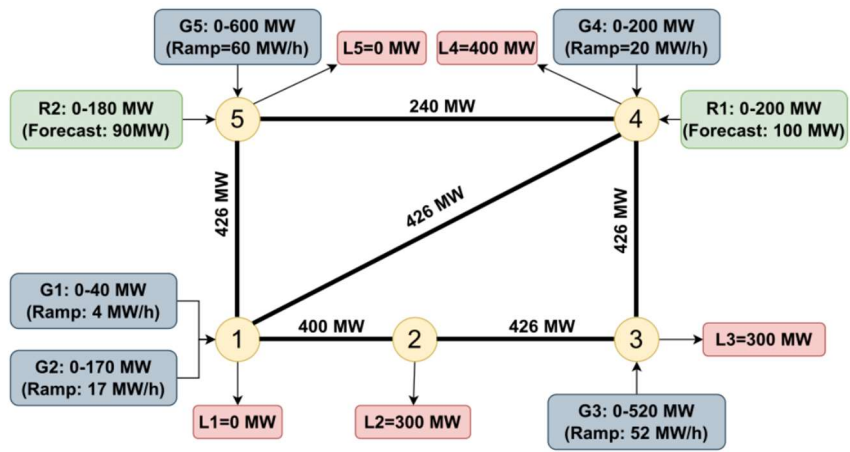

*Figure 2: The modified PJM 5-bus system with two renewable farms, R1 and R2, with forecasted active power 100MW and 90 MW, respectively.*

Data for the IEEE-39 bus system can be found in [15] and is also available from [12] and the associated tools. In the system modified for this paper, five renewable farms are connected: to buses 1, 5, 10, 15 and 20, with the predicted generation output of (respectively) 343 MW, 290 MW, 282 MW, 432 MW and 550 MW. The modified IEEE-118 bus system includes 20 renewable farms. The Matpower data file for all three modified case study systems can be found in [16].

To determine initial dispatch points, i.e., $p_i^*$ ($\forall i$) in (1d), a *reserve factor* of each non-renewable farm is assumed. This is a number between 0 and 1, which reduces the maximum available active power, and represent uncertainties in renewable farm available capacity. A deterministic DC-OPF calculation is run for each case study to calculate an initial dispatch point. Note that alternative initial dispatch points could be generated by other approaches, such as stochastic or robust optimisation.

Moreover, the number of initial guesses in $S_1$ is set as 100 for the IBLP approach, while the number of initial guesses in $S_2$ changes for each iteration $k$.

## 4.2 Simulation results

Simulation results of the PJM-5 bus system on CPT-01 are presented in Fig.3a. The solver takes 0.32 seconds within 6 iterations for the MILP approach and takes 0.48 seconds and 1.06 seconds with, respectively, one thread (IBLP-ST) and four threads (IBLP-MT) within 3 iterations. For the identified DRR shown in Fig.3a, its boundaries are connected to the network's binding operational constraints as: "FG"- the initial upper bound of R1, "DE"-the initial upper bound of R2, "AF"-initial lower bound of R2, "EG"- ramping down capability of G5, "BC"/"AB"/"CD"-thermal limit of Line 45.

From the power system operation perspective, if the realised renewable generation (R1, R2) falls within the *green area* (calculated DRR) in Fig.3a, the system operator can always find a dispatch strategy to balance the supply-demand without violating any operational constraints. If the DRR is calculated based on the revised PSD model under the electricity market environment, where renewables only participate in the energy market while FCAS are provided by other flexible generators (see Section 2), it can be instead interpreted as follows: for any realised renewable generation (R1,R2) in the *green area*, the desired FCAS (say regulation up/down FCAS) is adequate and can be dispatched by the system operator to compensate for deviations of both R1 and R2 from their forecasted values without incurring any operational violations.

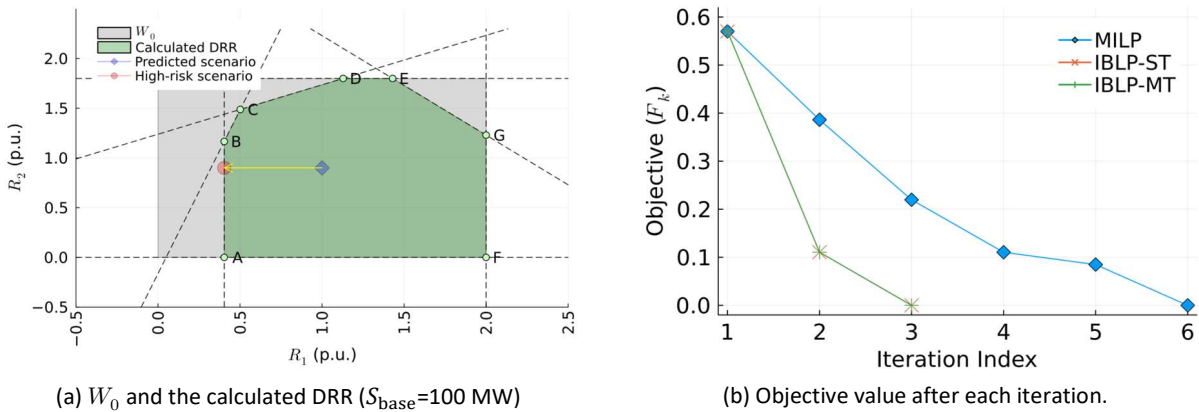

(a) $W_0$ and the calculated DRR ($S_{\text{base}}$=100 MW)   (b) Objective value after each iteration.

*Figure 3: The computational information in calculating DRR for the PJM-5 bus system (The yellow arrow represents the ramping* event *of highest risk in the left subfigure).*

The objective value of (4) in each iteration appears in Fig.3b, for the MILP, IBLP-ST, and IBLP-MT approaches. Using CPT-02, identical results are reported: both the number of iterations and the objective value after each iteration. The solution time, however, is 0.42, 0.30, and 1.79 seconds for the MILP, IBLP-ST and IBLP-MT approaches, respectively. The predicted renewable generation output, and the possible renewable output with the highest risk, appear in Fig.3a. The latter is identified as that falling on the DRR boundary and having the minimum distance to the *predicted scenario*. Correspondingly, the renewable generation ramping event with the highest risk and that should be paid most attention by the system operator, is the one where the renewable output moves from the *predicted* to the *high-risk* scenario.

*Table 1: Computational time for all studied cases with various approaches in both CPT-01 and CPT-02.*

| System Name | Platform | Computational Time (seconds) | | |
|---|---|---|---|---|
| | | MILP | IBLP-MT | IBLP-ST |
| PJM-5 Bus | CPT-01 | 0.32 | 0.48 | 1.05 |
| | CPT-02 | 0.42 | 0.30 | 1.79 |
| IEEE-39 Bus | CPT-01 | 30.69 | 7.66 | 15.57 |
| | CPT-02 | 56.47 | 7.08 | 27.04 |

The computational times for all studied cases with the various approaches in both CPT-01 and CPT-02 are also summarised in Table 1. For the IEEE-39 bus system, the objective values after each iteration for various approaches are presented in Fig.4a, where the MILP approach takes as high as 49 iterations while the IBLP approach only takes 10 iterations to report the solution. Although all approaches report the same DRR, the number of iterations taken by the MILP and IBLP approaches are significantly different.

Although the MILP approach guarantees a globally optimal DRR, these experiments suggest that it typically takes more iterations than the IBLP approach. The objective value reduces to zero slowly, i.e., only after several iterations. This could lead to unacceptable computational performance for large-scale systems with more renewable farms. That is, it may be computationally intractable or otherwise take an unacceptably long time to solve. In contrast, the number of iterations for the IBLP approach can be significantly lower since more than one linear cut can be generated to update $W$ in each iteration. Although this can lead to faster computation, the IBLP approach does not guarantee a global optimum. Nevertheless, using PCT in the IBLP approach can further improve computational speed, compare results between IBLP-ST and IBLP-MT. For the IEEE-39 bus system, the multi-threaded approach reduces computation time by 50.8% and 73.8% on CPT-01 and CPT-02, respectively.

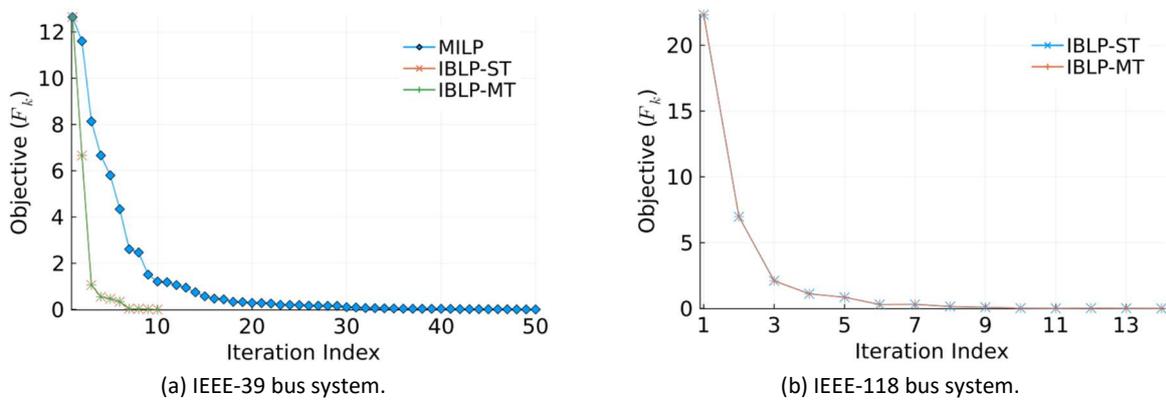

*Figure 4: Objective value after each iteration in calculating DRR for the IEEE-39 and IEEE-118 bus systems.*

Note further that although CPT-02 has both more CPU cores and a larger RAM than CPT-01, it takes more time to report solutions for both the MILP and IBLP-ST approaches. This is because the RAM on both CPT-01 and CPT-02 is sufficient, only one CPU is used for each of the MILP and IBLP-ST approaches, and the computing power of a single CPU core is greater for CPT-01 than for CPT-02.

The IEEE-118 bus system is further tested on CPT-01[6]. For this system with numerous buses and higher number of renewable farms, the MILP approach fails due to the computational requirements for solving a large MILP problem in each iteration. In contrast, the IBLP approaches calculate the DRR within only 14 iterations. The objective value after each iteration is presented in Fig.4b. The IBLP-ST and IBLP-MT take 394.21 seconds and 257.74 seconds, respectively, again demonstrating the benefit of PCT.

## 5. Conclusions and Potential Extensions

This paper studied the extent to which advanced computation techniques, including both HPC and PCT, could accelerate the computational speed in calculating DRR in power system operation. Three systems were studied on two computation platforms, testing two typical approaches, i.e., the MILP approach and the IBLP approach (including both IBLP-ST and IBLP-MT), showing that

1. Although the MILP approach guarantees a global optimum and is efficient for calculating DRR for small-scale systems with a few renewable farms, it can be computationally intractable for large-scale systems with a large number of renewable farms. In contrast, the IBLP approach appears more likely to solve within a reasonable time period, although it is not guaranteed to be globally optimal.

---

[6] As the remote server CPT-02 is shared with more broadly with other users, running Julia on CPT-02 encounters stability issues. Quite different computational times may result in otherwise identical runs. Thus, the IEEE-118 bus system is only tested on CPT-01.

2. The MILP approach is computationally challenging for large-scale systems with a high number of renewable farms (IEEE-118 bus system with 20 renewable farms in our case), even with HPC.
3. Applying PCT can significantly improve the computational performance of the IBLP approach. However, more experiments are required on real networks, which have an even larger number of buses and renewable farms.

Remarks and potential extensions of the DRR are summarised as follows.

1. Network reduction methods can produce an approximate power system model of smaller size. This may be useful if the *size* of a given problem is too large for *existing computational capability* to solve within a reasonable time. For example, several renewable farms can be approximated by a single renewable farm if they are sufficiently electrically close. Such reduction methods have been widely studied in both transmission and distribution networks [17–22].
2. Although in this paper, DRR was studied for renewable farms only, an extended dispatchable region that also considers loads at the transmission-distribution interface can clearly also be calculated. Such an extended dispatchable region provides information on how DERs in distribution networks should be constrained via, for example, operating envelopes [23–27].
3. As discussed in Section 2, DRR can also be applied within an electricity market environment, provided the ability to securely and reliably balance supply-demand can be appropriately quantified. High-risk ramping events, along with the corresponding binding operational constraints, can also be identified in this context.
4. A potential future research direction could be investigating DRR and how it could be interpreted if FCAS, particularly regulation down FCAS, is also available from renewable farms.

## 6. Appendix: Deriving the Equivalent Formulation of (4)

We first derive the dual problem of the inner *min* problem in (4), whose Lagrangian function can be expressed as follows [28].

$$
\begin{aligned}
L &= 1^T s + 1^T(t^+ + t^-) + \alpha^T(b - Aw - Bp - Cl - Dv - t^+ + t^-) + \delta^T(Ew + Fp + Gl \\
&\quad + Jv - s - d) \\
&= \alpha^T(b - Aw) + \delta^T(Ew - d) + (1 - \delta)^T s + (1 - \alpha)^T t^+ + (1 + \alpha)^T t^-
\end{aligned}
\tag{10}
$$

So, an equivalent formulation of the inner *min* problem in (4) is

$$\max_{\alpha,\delta} \min_{p,l,v} \alpha^T(b - Aw) + \delta^T(Ew - d) + (1 - \delta)^T s + (1 - \alpha)^T t^+ + (1 + \alpha)^T t^- \\
+ (-\alpha^T B + \delta^T F)p + (-\alpha^T C + \delta^T G)l + (-\alpha^T D + \delta^T J)v \tag{11a}$$

$$s.t. \quad s \geq 0, t^- \geq 0, t^+ \geq 0, \delta \geq 0 \tag{11b}$$

which can be further re-expressed as

$$\max_{\alpha,\delta} \alpha^T b - \delta^T d + (-\alpha^T A + \delta^T E)w \tag{12a}$$

$$s.t. \quad 0 \leq \delta \leq 1, -1 \leq \alpha \leq 1 \tag{12b}$$

$$\alpha^T B - \delta^T F = 0, \alpha^T C - \delta^T G = 0, \alpha^T D - \delta^T J = 0 \tag{12c}$$

Replacing the inner *min* problem in (4) by (12) leads to

$$\max_{\alpha,\delta} \alpha^T b - \delta^T d + \max_w (-\alpha^T A + \delta^T E)w \tag{13a}$$

$$s.t. \quad (5b) - (5c) \tag{13b}$$

$$Hw \leq f \tag{13c}$$